\documentclass[a4paper,12pt]{article}
\usepackage{amsthm}
\usepackage{amssymb}
\usepackage{amsmath}
\usepackage{amsfonts}
\usepackage{lineno}
\usepackage{color}

\newtheorem{remark}{Remark}[section]

\newtheorem{lemma}{Lemma}[section]
\newtheorem{theorem}{Theorem}[section]

\def\b1{\mbox{\boldmath $1$}}

\newenvironment{demo*}{\vspace{3mm}\noindent{\bf Proof.}}{\hfill $\Box$ \vspace{3mm}}

\begin{document}

\title{\bf \Large Finite-time ruin probabilities of bidimensional risk models with correlated Brownian motions  }
\author{\normalsize {Dan Zhu $^{1}$},\ {Ming Zhou $^{2}$} and {Chuancun Yin$^{1,*}$}\\
\normalsize $^{1}$ School of Statistics and Data Science, Qufu Normal University,\\
\noindent{\normalsize Shandong 273165, China; zhudanspring@qfnu.edu.cn (D.Z.)}\\
\normalsize $^{2}$ Center for Applied Statistics, School of Statistics, Renmin University of China, \\
\noindent{\normalsize Beijing 100872, China; mingzhou@ruc.edu.cn (M.Z.)}\\
\normalsize {\ Correspondence: ccyin@qfnu.edu.cn (C.Y.)}}

 \maketitle
 \noindent{\large {\bf Abstract}}
 \textcolor{blue} {The present work concerns the finite-time ruin probabilities for several bidimensional risk models with constant interest force and correlated Brownian motions.} Under the condition that the two Brownian motions $\{B_1(t), t\ge 0\}$ and  $\{B_2(t), t\ge 0\}$  are correlated, we  establish new results for the finite-time ruin probabilities.  \textcolor{blue} {Our research has enriched the development of the ruin theory with heavy tails in unidimensional risk models and the dependence theory of stochastic processes.}

\medskip


\newpage

\baselineskip=20pt

\section{Introduction}\label{intro}

\textcolor{blue}{ In traditional studies, many researchers investigated the ruin probabilities problems of an insurer under unidimensional models. For example,  Tang (2004) studied the ruin probability problems with constant interest force. Other studies about these problems can also be seen in Tang (2005), Hao and Tang (2008), Wang (2008), Wang et al. (2012). An assumption behind this is that the
businesses of the insurer are homogeneous and can be described by a unidimensional model. The assumption is too strong. Thus, bidimensional or multidimensional insurance risk models have received growing interest in recent years, such as Li (2017), Cheng (2019) and Chen et al. (2019). Various assumptions have been considered regarding the claim arrival processes and the claim amount distributions. See e.g. Chan et al. (2003), Li et al. (2007), Avram et al. (2008) and Chen et al. (2011). Chen et al. (2013) considered finite-time ruin probabilities of nonstandard bidimensional renewal risk models with constant interest
forces and diffusion generated by Brownian motions, they assumed that the two Brownian motions $\{B_1(t), t\ge 0\}$ and  $\{B_2(t), t\ge 0\}$ are mutually independent.  Similar results were obtained by Yang and Li (2014), but they consider dependent subexponential claims.
 More papers can be seen in Wang et al. (2022), Lu and Yuan (2022) and the references therein.
In this paper, we consider
uniform asymptotics for the finite-time ruin probabilities for several bidimensional risks models with constant interest force and correlated Brownian motions that means that the businesses of the insurer have a relationship with each other. We will introduce the risk models and different types of ruin times corresponding ruin probabilities as follows.}

 The bidimensional risk model $\vec {U}(t)=(U_1(t),U_2(t))^{\tau}$ is the surplus vector of an insurance company at time $t\geq0$, which we consider in this paper can be formally stated as
\begin{equation}
U_i(t)=u_i e^{rt}+\int_0^t e^{r(t-s)}dC_i(s)-\int_0^t e^{r(t-s)}dS_i(s)+\sigma_i \int_0^t e^{r(t-s)}dB_i(s),\ t\ge 0,
\end{equation}
where $\vec {u}=(u_1,u_2)^{\tau}$ stands for the initial surplus vector, $\vec {C}(t)=(C_1(t),C_2(t))^{\tau}$ for the total premiums received up to time $t$, and $\{C_1(t), t\ge 0\}$, $\{C_2(t), t\ge 0\}$ are mutually independent, $r\ge 0$ for the interest rate, $(S_1(t), S_2(t))=(\sum_{i=1}^{N_1(t)}X_{1i}, \sum_{i=1}^{N_2(t)}X_{2i})$ for the total amount of claims vector up to time $t$. Here
$\vec {X_i}=(X_{1i}, X_{2i})^{\tau}, i=1,2,\cdots$, denote pairs of claims whose  arrival times constitute a counting process vector $\{{\vec{N}}(t), t\ge 0\}$, where ${\vec{N}}(t)=(N_1(t), N_2(t))$, and $\{N_1(t), t\ge 0\}$, $\{N_2(t), t\ge 0\}$ are mutually independent. The process $\{N_i(t), t\ge 0\}$ is a Poisson process with intensity $\lambda_i>0$ and $\{\vec{X_i}, i=1,2,\cdots\}$ is a sequence of independent copies of the random pair $\vec {X}=(X_1,X_2)^{\tau}$ with joint distribution function $F(x_1, x_2)$ and marginal distribution functions $F_1(x_1)$ and $F_2(x_2)$. All vectors $\vec {X_i}$'s and $\vec {C}$ consist of only nonnegative components,   $\vec {C}(0)=(0,0)^{\tau}$. Moreover, each $C_i(t)$ is a nondecreasing and right-continuous stochastic process. The vector  $\vec {B}(t)=(B_1(t),B_2(t))^{\tau}$ denotes a standard bidimensional Brownian motion with constant correlation coefficient $\rho\in [-1,1]$, while $\sigma_1\ge 0$ and $\sigma_2\ge 0$ are constants.  For simplicity, we assume that $\{\vec {X_i}, i=1,2,\cdots\}$,  $\{{\vec{N}}(t), t\ge 0\}$ and $\{\vec {C}(t), t\ge 0\}$  are independent, and furthermore, both of them are also independent of $\{\vec {B}(t), t\ge 0\}$.
To avoid the certainty of ruin in each class, we assume that the following safety
loading conditions hold when $r=0$:
$$EC_i(t)-\lambda_i EX_{i1}>0, \ i=1, 2.$$

 In this paper, we consider the following four types of ruin probabilities:
For finite-horizon $T>0$, we  define
 \begin{equation}
 \psi_{\max}(\vec {u},T)=P(T_{\max}\le T|\vec{U}(0)=\vec{u}),
 \end{equation}
 where
 $$T_{\max}=\inf\{t>0|\max\{U_1(t),U_2(t)\}<0\};$$
\begin{equation}
 \psi_{\min}(\vec {u},T)=P(T_{\min}\le T|\vec{U}(0)=\vec{u}),
 \end{equation}
 where
 $$T_{\min}=\inf\{t>0|\min\{U_1(t),U_2(t)\}<0\};$$
 and
\begin{equation}
 \psi_{\rm{sum}}(\vec {u},T)=P(T_{\rm {sum}}\le T|\vec{U}(0)=\vec{u}),
 \end{equation}
 where
 $$T_{\rm{sum}}=\inf\{t>0|U_1(t)+U_2(t)<0\};$$
 \begin{equation}
 \psi_{\rm and}(\vec {u},T)=P(T_{\rm {and}}\le T|\vec{U}(0)=\vec{u}),
 \end{equation}
 where $T_{\rm {and}}=\max\{T_1, T_2\}$ and
 $$T_{i}=\inf\{t>0|U_i(t)<0 \ {\rm for \ some} \ 0\le t\le T),\ i=1,2,$$
with $\inf\emptyset=\infty$ by convention.

We remark that the  probability in (1.2) denotes the probability of ruin  \textcolor{blue} {occurs when} both $U_1(t)$ and $U_2(t)$ are below zero at the same time  within finite time  $T>0$, the  probability in (1.3) denotes the probability of ruin occurs that at least one of  $\{U_i(t), i=1,2\}$ is below zero   within finite time  $T>0$, the  probability in (1.4) denotes the probability of ruin occurs that the total of  $U_1(t)$ and $U_2(t)$ is negative   within finite time  $T>0$, whereas the  probability in (1.5)
denotes the probability of ruin occurs that both $U_1(t)$ and $U_2(t)$ are below zero, not necessarily  simultaneously, within finite time  $T>0$.
$T_{\rm {and}}$ represents a more critical time than $T_{\max}$,
  and the ruin probability defined by $T_{\rm sum}$ will be reduced to that in the unidimensional model. The following relation holds between the four  ruin probabilities defined above:
 $$\psi_{\max}(\vec {u},T)\le \psi_{\rm and}(\vec {u},T)\le\psi_{\min}(\vec {u},T),\ \psi_{\rm{sum}}(\vec {u},T)\le \psi_{\min}(\vec {u},T),$$
 and
\begin{equation} \psi_{\min}(\vec {u},T)+\psi_{\rm and}(\vec {u},T)=P(T_1\le T|U_1(0)=u_1)
+P(T_2\le T|U_2(0)=u_2).
\end{equation}
\textcolor{blue} {The rest of this paper is organized as follows:} In Section 2, we review the related results after briefly introducing some preliminaries about heavy-tailed distributions; in Section 3, we give some definitions and lemmas;
the main results and the proofs procedure are given in Section 4.

 \vskip 0.2cm
\section{ Review of related results}
\setcounter{equation}{0}

 \textcolor{blue} {Unless otherwise stated herein, all limits relations are for $(u_1, u_2)\to (\infty,\infty)$. Denote $a\lesssim b$ and  $a \gtrsim b$ if
$\limsup a/b \le 1$ and $\limsup a/b \ge 1$,respectively,$a\sim b$ if both, where, $a(\cdot, \cdot)$ and  $b(\cdot, \cdot)$ are two positive functions. Let  $F_1*\cdots*F_n$ as the convolution of the distributions $F_1, \cdots, F_n$ and by $F^{*n}$ the $n$-fold convolution of a distribution $F$.}

\textcolor{blue}{In this section, we will review some definitions and properties that are relevant to the results of this paper and consider only the case of the distribution of heavy-tail claims.}
A r.v. $X$ or its d.f. $F(x)=1-\overline{F}(x)$ satisfying
$\overline{F}(x)>0$ for all $x\in (-\infty,\infty)$ is called
heavy-tailed to the right, or simply heavy-tailed, if $E[e^{\gamma
X}]=\infty$ for all $\gamma>0$. We recall here some important
classes of heavy-tailed distributions as follows.

\textcolor{blue}{ $F$ is a long tailed distribution, written as $F\in\mathcal{L}$, if for some $t>0$,}
$\lim\limits_{x\rightarrow\infty}\frac{\overline{F}(x-t)}{\overline{F}(x)}=1 $
holds. Note  that the convergence  is  uniform over
$t$ in compact intervals.
\textcolor{blue}{If
$\lim\limits_{x\rightarrow\infty}\frac{\overline{F^{\ast n}}(x)}{\overline{F}(x)}=n$
 holds $(n=2, 3,\cdots)$, $F$ is a subexponential distribution on $(0,\infty)$,
written as $F\in \mathcal{S}$.
 For some  $0<t<1$, if
$\limsup\limits_{x\rightarrow\infty}\frac{\overline{F}(tx)}{\overline{F}(x)}<\infty $
holds, $F$ is said to be dominatedly varying tailed distribution, written as $F\in
\mathcal{D}$. We call
$F$ a consistently varying tailed distribution, }written as $F\in
\mathcal{C}$, if
$$\lim_{t\downarrow1}\liminf_{x\rightarrow\infty}\frac{\overline{F}(tx)}{\overline{F}(x)}=1,\ \ {\text{\rm
 or, equivalently}}, \ \
\lim_{t\uparrow1}\limsup_{x\rightarrow\infty}\frac{\overline{F}(tx)}{\overline{F}(x)}=1$$
holds. A distribution
$F$ is a extended regularly varying tailed, written as
$F\in \mathcal{ERV}(-\alpha, -\beta)$, for some $0\le\alpha\leq\beta<\infty$, if for $s\geq1$,
$s^{-\beta}\leq\liminf\limits_{x\rightarrow\infty}\frac{\overline{F}(sx)}{\overline{F}(x)}\leq
\limsup\limits_{x\rightarrow\infty}\frac{\overline{F}(sx)}{\overline{F}(x)}\leq
s^{-\alpha} $holds.

\textcolor{blue} {It is obvious that the following formula holds}
$$ \mathcal{ERV}(-\alpha, -\beta)\subset \mathcal{C}\subset
\mathcal{D}\cap \mathcal{L}\subset \mathcal{S}\subset
\mathcal{L}. $$
\textcolor{blue}{There are many other references to heavy-tailed distributions, the reader can also refer to} Bingham et al.
(1987), Cline and Samorodnitsky (1994),  Tang and Tsitsiashvili
(2003), Asmussen and Albrecher (2010), Konstantin (2018), \textcolor{blue}{Ke et al. (2019)}, among others.

The asymptotic behavior of the finite-time ruin probability of bidimensional or multidimensional risk models has been investigated by
Liu et al. (2007). They proved that under the conditions $F_1, F_2\in \mathcal{S}$, $N_1(t)=N_2(t)$, $\sigma_1=\sigma_2=0$, $r>0$ and  the claim vector $\vec{X}$  consist of independent components,
$$ \psi_{\max}(\vec{u};T)\sim\frac{\lambda(\lambda+\frac{1}{T})}{r^2}\int_{u_1}^{u_1 e^{rT}}\frac{\overline{F_1}(y)}{y}dy
\int_{u_2}^{u_2 e^{rT}}\frac{\overline{F_2}(y)}{y}dy,\ \  \text{\rm as}\ \ (u_1,u_2)\to (\infty,\infty).$$
Under the conditions $F_1, F_2\in \mathcal{S}$, $r=0$, $N_1(t)=N_2(t)$, $C_i(\cdot)$ are deterministic linear functions   and both the claim vector $\vec{X}$ and the bidimensional Brownian motion $\vec{B}$ consist of independent components, Li, Liu and Tang (2007) obtained that, for each fixed time $T>0$,
\begin{equation}
 \psi_{\max}(\vec{u};T)\sim \lambda T(1+\lambda T)\overline{F_1}(u_1)\overline{F_2}(u_2),\ \  \text{\rm as}\ \ (u_1,u_2)\to (\infty,\infty).\nonumber
 \end{equation}
 Chen, Yuen and Ng (2011) investigated the uniform asymptotics of  $\psi_{\rm and}(\vec {u},T)$ and $\psi_{\min}(\vec {u},T)$ for  an ordinary renewal risk model  with the claim amounts belonging to the consistently varying
tailed distributions class for large $T$. Zhang and Wang (2012) considered model (1.1) with $r=0$ and assumed that
all sources of randomness, $\{X_{1k}, k=1,2,\cdots\}$, $\{X_{2k}, k=1,2,\cdots\}$, $\{N_1(t)=N_2(t), t\ge 0\}$, $\{B_1(t), t\ge 0\}$ and  $\{B_2(t), t\ge 0\}$
are mutually independent. They obtained that if $F_1, F_2\in \mathcal{EVR}(-\alpha,-\beta)$ for some $0<\alpha\le \beta<\infty$, then, for each fixed time $T\ge 0$,
\begin{equation}
 \psi_{\max}(\vec{u};T)\sim \lambda T(1+\lambda T)\overline{F_1}(u_1)\overline{F_2}(u_2),\ \  \text{\rm as}\ \ (u_1,u_2)\to (\infty,\infty).\nonumber
 \end{equation}
 The analogous result for  multidimensional  risk model can be found in Asmussen and Albrecher (2010, P.441).
\vskip 0.2cm
\section{Some lemmas}
\setcounter{equation}{0}
Before giving the main results, we first provide some lemmas.

\begin{lemma} If $F\in\mathcal{S}$, then for each $\varepsilon>0$, there exists some constant $C_{\varepsilon}>0$ such that the inequality
$$\overline{F^{*n}}(x)\le  C_{\varepsilon}(1+\varepsilon)^n \overline{F}(x)$$
 holds for all $n=1,2,\cdots$ and $x\ge 0$.
\end{lemma}
{\bf Proof}. See Lemma 1.3.5 of Embrechts et al. (1997).
\begin{lemma}  Let $G_1$ and $G_2$ be two distribution functions. If $G_1\in\mathcal{S}$ and $\overline{G_2}(x)=o(\overline{G_1}(x))$,
 then we have  $\overline{G_1*G_2}(x)\sim \overline{G_1}(x)$ as $x\to\infty$.
\end{lemma}
{\bf Proof}. See Proposition 1 of Embrechts et al. (1997).
\begin{lemma} Consider a unidimensional  risk model
\begin{equation}
U_i(t)=u_i + C_i(t)-S_i(t)+\sigma_i B_i(t),\ t\ge 0, i=1,2.
\end{equation}
 If $F_i\in\mathcal{S}$,
then the ruin probability with finite-horizon $T$ satisfies that
$$\psi_i(u_i; T)=P(U_i(t)<0 \ {\rm for}\ {\rm  some}\  t\le T|U_i(0)=u_i)\sim  \lambda T\overline{F_i}(u_i), u_i\to\infty.$$
\end{lemma}
{\bf Proof}. Clearly, on the one hand,
\begin{eqnarray}
 \psi_i(u_i; T)&\ge& P(S_i(T)\ge u_i+C_i(T)+\sigma_i \sup_{0\le t\le T}B_i(t))\nonumber\\
 &=&\int_0^{\infty}P(S_i(T)\ge u_i+C_i(T)+\sigma_i z)dP(\sup_{0\le t\le T}B_i(t))\le z)\nonumber\\
 &=&P(S_i(T)\ge u_i)\int_0^{\infty}\int_0^{\infty}\frac{P(S_i(T)\ge u_i+l_i+\sigma_i z)}{P(S_i(T)\ge u_i)}dP(\sup_{0\le t\le T}B_i(t))\le z)\nonumber\\
&& \times dP(C_i(T)\le  l_i)\nonumber\\
 &\sim& P(S_i(T)\ge u_i),
 \end{eqnarray}
where we have used the fact that $P(S_i(T)\ge u_i+l_i+\sigma_i z)\le P(S_i(T)\ge u_i)$ and the dominated convergence theorem.

On the other hand,
\begin{eqnarray}
 \psi_i(u_i; T)&\le& P(S_i(T)+\sigma_i \sup_{0\le t\le T}(-B_i(t))\ge u_i)\nonumber\\
 &\sim& P(S_i(T)\ge u_i),
 \end{eqnarray}
where we have used Lemma 3.2 and the fact that
$$ P(\sigma_i \sup_{0\le t\le T}(-B_i(t))\ge u_i)=o(P(S_i(T)\ge u_i)).$$
By Lemma 3.1 and dominated convergence theorem, we have
$$ P(S_i(T)\ge u_i)\sim \overline{F_i}(u_i) \sum_{n=1}^{\infty}nP(N(T)=n)=\lambda T  \overline{F_i}(u_i), \ as\  u_i\to\infty.$$
The result follows from (3.2) and (3.3).

\begin{lemma} Consider a unidimensional  risk model
\begin{equation}
 U_i(t)=u_i e^{rt}+\int_0^t e^{r(t-s)}C_i(ds)-\int_0^t e^{r(t-s)}dS_i(s)+\sigma_i \int_0^t e^{r(t-s)}dB_i(s),\ t\ge 0, i=1,2.\nonumber
\end{equation}
 If $F_i\in\mathcal{S}$,
then the ruin probability with finite-horizon $T$ satisfies that
$$\psi_i(u_i; T)=P(U_i(t)<0 \ {\rm for}\ {\rm  some}\  t\le T|U_i(0)=u_i)\sim
\frac{\lambda}{r}\int_{u_i}^{u_i e^{rT}}\frac{\overline{F_i}(y)}{y}dy,\  u_i\to\infty.$$
\end{lemma}
{\bf Proof}. Just modify the proof of Lemma 3.3 we have
$$\psi_i(u_i; T)\sim P\left(\sum_{j=1}^{N(T)} X_{ij}e^{-r\tau_j}\ge u_i\right)\sim \lambda \int_0^{T}P(X_{i1}e^{-rz}>u_i)dz,\  u_i\to\infty,$$
 in the last step, we used (4.14) in Tang (2005). Here $\tau_j$ are the arrival times of the Poisson process $N(t)$.
 \textcolor{blue} {In fact, $$z=\frac{1}{r}\log\frac{y}{u_i},$$ we have that
$$dz=d\left(\frac{1}{r}\log\frac{y}{u_i}\right)=\frac{1}{r}\cdot\frac{u_i}{y}\cdot\frac{1}{u_i}dy=\frac{1}{ry}dy,$$
Then, $$\lambda\int_0^TP(X_{i1}>u_ie^{rz})dz=\frac{\lambda}{r}\int_{u_i}^{u_ie^{rT}}\frac{\overline{F}_i(y)}{y}dy.$$}

Upon a trivial substitution implies the required result.

{\bf Defination 3.1} (i)\ Two processes $\{X_1(t); t\ge 0\}$ and $\{X_2(t); t\ge 0\}$ are said to be  positively  associated  if
$${\rm Cov}(f(X_1(t_1),X_2(t_2)), g(X_1(t_1),X_2(t_2))|X_1(0)=x_1, X_2(0)=x_2)\ge  0,$$
for all non-decreasing real valued functions $f$ and $g$ such that the covariance exists, all $t_1, t_2\ge 0$, and all $x_1, x_2\in \Bbb{R}$.

(ii)\ Two processes $\{X_1(t); t\ge 0\}$ and $\{X_2(t); t\ge 0\}$ are said to be  negatively  associated  if
$${\rm Cov}(f(X_1(t_1)), g(X_2(t_2))|X_1(0)=x_1, X_2(0)=x_2)\le  0,$$
for all non-decreasing real valued functions $f$ and $g$ such that the covariance exists, all $t_1, t_2\ge 0$, and all $x_1, x_2\in \Bbb{R}$.

 {\bf Defination 3.2} Two processes $\{X_1(t); t\ge 0\}$ and $\{X_2(t); t\ge 0\}$ are said to be  positively (negatively) quadrant dependent  if
\begin{eqnarray}
&& P(X_1(t_1)>y_1,X_2(t_2)>y_2|X_1(0)=x_1, X_2(0)=x_2)\nonumber\\
 &&\ge (\le)  P(X_1(t_1)>y_1 |X_1(0)=x_1) P(X_2(t_2)>y_2 |X_2(0)=x_2),
\end{eqnarray}
for all $t_1, t_2\ge 0$, and for all $y_1, y_2, x_1, x_2\in \Bbb{R}$.

It is well known that (cf. Ebrahimi (2002))  $(X_1(t), X_2(t))$ is positively (negatively)  associated implies
 $X_1(t)$ and $X_2(t)$ are    positively (negatively) quadrant dependent.

 Let  $\vec {B}(t)=(B_1(t),B_2(t))^{\tau}$ be a standard bidimensional Brownian motion with constant correlation coefficient $\rho\in (-1,1)$.
 For notional convenience, for $t\ge 0$, write
$\underline{B}_i(t)=\inf_{0\le s\le t}B_i(s),\ \overline{B}_i(t)=\sup_{0\le s\le t}B_i(s),\ i=1,2.$ It is well known that, for $x>0$,
$P(\underline{B}_i(t)<-x)=P(\overline{B}_i(t)>x)=2P(B_i(t)>x).$
The following lemma is essential to prove our main results. Moreover, it is of independent interest.
\begin{lemma} For any $x_1>0, x_2>0$, if  $\rho\in [0, 1)$, then
\begin{equation}
P(\overline{B}_1(t)>x_1, \overline{B}_2(t)>x_2)\ge P(\overline{B}_1(t)>x_1) P(\overline{B}_2(t)>x_2),
\end{equation}
and
\begin{equation}
P(\underline{B}_1(t)<-x_1, \underline{B}_2(t)<-x_2)\ge P(\underline{B}_1(t)<-x_1) P(\underline{B}_2(t)<-x_2);
\end{equation}
If $\rho\in (-1,0]$, then
\begin{equation}P(\overline{B}_1(t)>x_1, \overline{B}_2(t)>x_2)\le P(\overline{B}_1(t)>x_1) P(\overline{B}_2(t)>x_2)
 \end{equation}
and
\begin{equation}
P(\underline{B}_1(t)<-x_1, \underline{B}_2(t)<-x_2)\le P(\underline{B}_1(t)<-x_1) P(\underline{B}_2(t)<-x_2).
\end{equation}
\end{lemma}

{\bf Proof}. For any  $t_1, t_2\ge 0$, we have ${\rm Cov}(B_1(t_1), B_2(t_2))=\rho\min\{t_1, t_2\}$.
It follows from the Theorem in Pitt (1982) that $\rho\ge 0$ is necessary and sufficient for $(B_1(t),B_2(t))^{\tau}$ to be positively associated since $(B_1(t_1), B_2(t_2))^{\tau}$ is bivariate normal, which implies  that $(B_1(t),B_2(t))^{\tau}$  is positively  quadrant dependent.
Thus (3.5) holds.
To prove (3.6), we use (3.5) and the  facts that
$-\sup_{0\le s\le t}B_i(s)=\inf_{0\le s\le t}(-B_i(s))$ and $(-B_1(t), -B_2(t))^{\tau}$ is also a standard bidimensional Brownian motion with correlation coefficient $\rho$. Inequalities (3.7) and (3.8) can be proved similarly. This completes the proof.

  For \textcolor{blue} { $r\geq0$, }consider a bidimensional Gaussian process  $(\int_0^t e^{-rs}dB_1(s),\int_0^t e^{-rs}d B_2(s))^{\tau}$, where $\vec {B}(t)=(B_1(t),B_2(t))^{\tau}$ is a standard bidimensional Brownian motion with constant correlation coefficient $\rho\in (-1,1)$.
 For $t\ge 0$, write
\textcolor{blue} { $$\underline{\Delta}_i(t)=\inf_{0\le s\le t} \int_0^s e^{-rl}dB_1(l),\ \overline{\Delta}_i(t)=\sup_{0\le s\le t}\int_0^s e^{-rl}d B_2(l)
,\ i=1,2.$$}
The following lemma is an extension of Lemma 3.5.
\begin{lemma} For any $x_1>0, x_2>0$, if  $\rho\in [0, 1)$, then
 $$P\left(\overline{\Delta}_1(t)>x_1,   \overline{\Delta}_2(t)>x_2\right)
\ge P\left(\overline{\Delta}_1(t) >x_1\right) P\left(\overline{\Delta}_2(t)>x_2\right), $$
and
 $$
P(\underline{\Delta}_1(t)<-x_1,  \underline{\Delta}_2(t)<-x_2)\ge P(\underline{\Delta}_1(t)<-x_1) P(\underline{\Delta}_2(t)<-x_2);
 $$
If $\rho\in (-1,0]$, then
  $$P\left(\overline{\Delta}_1(t)>x_1,   \overline{\Delta}_2(t)>x_2\right)
\le P\left(\overline{\Delta}_1(t) >x_1\right) P\left(\overline{\Delta}_2(t)>x_2\right), $$
and
 $$
P(\underline{\Delta}_1(t)<-x_1,  \underline{\Delta}_2(t)<-x_2)\le P(\underline{\Delta}_1(t)<-x_1) P(\underline{\Delta}_2(t)<-x_2).
 $$
\end{lemma}

\begin{remark}
Several distributions of interest are available in closed form (see, e.g. He, Keirstead and  Rebholz (1998)).
These include the joint distributions of  $(\underline{X}_1(t), \underline{X}_2(t))$,  $(\overline{X}_1(t), \overline{X}_2(t))$,
$(\underline{X}_1(t), \overline{X}_1(t))$, and so on. But those closed-form results can not apply our proofs to the main results. The results of Lemmas 3.5 and 3.6 can not be obtained to the results of Shao and Wang (2013).
\end{remark}

\begin{lemma} Let $\{N(t), t\ge 0\}$ be a Poisson process with arrival times $\tau_k, k=1,2,\cdots$. Given $N(T)=n$ for arbitrarily fixed $T>0$ and $n=1,2,\cdots$, the random vector $(\tau_1,\cdots,\tau_n)$ is equal in distribution to the random vector $(TU_{(1,n)}, \cdots, TU_{(n,n)})$, where  $U_{(1,n)},\cdots, U_{(n,n)}$ denote the order statistics of $n$  i.i.d. (0,1) uniformly distributed random variables $U_{1},\cdots, U_{n}$.
\end{lemma}
{\bf Proof.} See Theorem 2.3.1 of Ross (1983).

\begin{lemma} Let $X$ and $Y$ be two independent and nonnegative random variables. If $X$ is subexponentially distributed while $Y$ is bounded and nondegenerate at $0$, then the product $XY$ is subexponentially distributed.
\end{lemma}
{\bf Proof.} See Corollary 2.3 of Cline and Samorodnitsky (1994).

The following result is due to Tang (2004).
\begin{lemma} Let $X$ and $Y$ be two independent random variables with distributions $F_X$ and $F_Y$. Moreover, $Y$ is  nonnegative and nondegenerate at $0$. Then
$$F_{X-Y}\in \mathcal{L}\Leftrightarrow F_X\in \mathcal{L}\Leftrightarrow \overline{F}_{X-Y}(x)\sim \overline{F}_{X}(x).$$
\end{lemma}

\section{Main results and proofs}
\setcounter{equation}{0}
In this paper, we establish new results for the finite-time ruin probabilities. Unlike the above-motioned articles, we assume that the two Brownian motions $\{B_1(t), t\ge 0\}$ and  $\{B_2(t), t\ge 0\}$ are correlated with constant correlation coefficient $\rho\in (-1,1)$.
The followings are the main results of this paper.

\begin{theorem} Consider the insurance risk model introduced in Section 1. Assume that $N_1(t)=N_2(t)=N(t)$, $\rho\in (-1,0]$, $r=0$ and $\{X_{1k}, k=1,2,\cdots\}$, $\{X_{2k}, k=1,2,\cdots\}$, $\{C_1(t), t\ge 0\}$, $\{C_2(t), t\ge 0\}$,  $\{N(t), t\ge 0\}$, $\{(B_1(t), B_2(t)), t\ge 0\}$
are mutually independent.\\
(a)\ If  $F_1, F_2\in \mathcal{S}$, then, for each fixed time $T\ge 0$,
\begin{equation}
 \psi_{\max}(\vec{u};T)\sim \lambda T(1+\lambda T)\overline{F_1}(u_1)\overline{F_2}(u_2),\ \  \text{\rm as}\ \ (u_1,u_2)\to (\infty,\infty),
 \end{equation}
 \begin{equation}
 \psi_{\min}(\vec{u};T)\sim \lambda T\left(\overline{F_1}(u_1)+\overline{F_2}(u_2)\right),\ \  \text{\rm as}\ \ (u_1,u_2)\to (\infty,\infty).
 \end{equation}
 (b)\ If  $F_1*F_2\in \mathcal{S}$, then, for each fixed time $T\ge 0$,
 \begin{equation}
 \psi_{\rm{sum}}(\vec{u};T)\sim \lambda T\left(\overline{F_1}(u_1+u_2)+\overline{F_2}(u_1+u_2)\right),\ \  \text{\rm as}\ \ u_1+u_2\to \infty.
 \end{equation}
\end{theorem}

{\bf Proof. }\ First, we establish the asymptotic upper bound for $\psi_{\max}(\vec{u};T)$.
Clearly,
 \begin{eqnarray}
 \psi_{\max}(\vec{u};T)&\le& P\left(\sum_{i=1}^{N(T)}\vec{X}_i-{\sigma_1 \underline{B}_1(T)\choose
  \sigma_2 \underline{B}_2(T)}>\vec{u}\right)\nonumber\\
 &=&\sum_{n=0}^{\infty}P(N(T)=n)P\left(\sum_{i=1}^{n}\vec{X}_i- {\sigma_1 \underline{B}_1(T)\choose
  \sigma_2 \underline{B}_2(T)}>\vec{u}\right)\nonumber\\
  &=&\sum_{n=0}^{\infty}P(N(T)=n)\int_0^{\infty}\int_0^{\infty}P(\sum_{i=1}^{n}\vec{X}_i\in d\vec{z})\nonumber\\
  &&\times P\left(\vec{z}-{\sigma_1 \underline{B}_1(T)\choose
  \sigma_2 \underline{B}_2(T)}>\vec{u}\right).
\end{eqnarray}
Since $\rho\in (-1,0]$, using (3.8) one has
\begin{eqnarray}
&&P\left(\vec{z}-{\sigma_1 \underline{B}_1(T)\choose
\sigma_2 \underline{B}_2(T)}>\vec{u}\right)\nonumber\\
&&\le P(z_1-\sigma_1 \underline{B}_1(T)>u_1)P(z_2-\sigma_2\underline{B}_2(T)>u_2).
\end{eqnarray}
Using the independence of $\{X_{1k}, k=1,2,\cdots\}$ and  $\{X_{2k}, k=1,2,\cdots\}$ one has
\begin{equation}
P\left(\sum_{i=1}^{n}\vec{X}_i\in\vec{z}\right)=P\left(\sum_{i=1}^{n}X_{1i}\in dz_1\right)P\left(\sum_{i=1}^{n}X_{2i}\in dz_2\right).
\end{equation}
Substituting (4.5) and (4.6) into (4.4) and using the  dominated convergence theorem, we get
 \begin{eqnarray}
\psi_{\max}(\vec{u};T)&\le& \sum_{n=0}^{\infty}P(N(T)=n) P\left(\sum_{i=1}^{n}X_{1i}-\sigma_1 \underline{B}_1(T)>u_1\right)
  P\left(\sum_{i=1}^{n}X_{2i}-\sigma_2 \underline{B}_2(T)>u_1\right)\nonumber\\
  &\sim& \sum_{n=0}^{\infty}P(N(T)=n)n^2\overline{F_1}(u_1)\overline{F_2}(u_2)\nonumber\\
 &=& \lambda T(1+\lambda T)\overline{F_1}(u_1)\overline{F_2}(u_2),\ \  \text{\rm as}\ \ (u_1,u_2)\to (\infty,\infty),
 \end{eqnarray}
where in the second step, we have used Lemma 3.2 and the fact that
$$ P\left(\sigma_j \sup_{0\le t\le T}(-B_j(t))\ge u_j\right)=o\left(P(\sum_{i=1}^n X_{ji}\ge u_j)\right), \ j=1,2.$$
Next,  we establish the asymptotic lower bound for $\psi_{\max}(\vec{u};T)$.
Clearly,
 \begin{eqnarray}
 \psi_{\max}(\vec{u};T)&\ge& P\left(\sum_{i=1}^{N(T)}\vec{X}_i-\vec{C}(T)-{\sigma_1 \overline{B}_1(T)\choose
  \sigma_2 \overline{B}_2(T)}>\vec{u}\right)\nonumber\\
 &=&\sum_{n=0}^{\infty}P(N(T)=n)P\left(\sum_{i=1}^{n}\vec{X}_i- {\sigma_1 \overline{B}_1(T)\choose
  \sigma_2 \overline{B}_2(T)}-\vec{C}(T)>\vec{u}\right)\nonumber\\
  &\equiv& \sum_{n=0}^{\infty}P(N(T)=n) I_1,
\end{eqnarray}
where $I_1$ can be written as
 \begin{equation}
I_1=\int_0^{\infty}\int_0^{\infty}P(\overline{B_1}(T)\in dy_1, \overline{B_2}(T)\in dy_2)J_1 J_2.
\end{equation}
Here
$$J_1=P\left(\sum_{i=1}^{n}X_{1i}-C_1(T)-\sigma_1 y_1>u_1\right),$$
and
$$J_2=P\left(\sum_{i=1}^{n}X_{2i}-C_2(T)-\sigma_2 y_2>u_2\right).$$
 For large constants  $a>0$ and $b>0$, we further write $I_1$ as
\begin{eqnarray}
 I_1&=&\left(\int_0^{a}\int_0^{b}+\int_0^{a}\int_b^{\infty}+\int_a^{\infty}\int_0^{b}+\int_a^{\infty}\int_b^{\infty}\right)P(\overline{B_1}(T)\in dy_1, \overline{B_2}(T)\in dy_2)J_1 J_2\nonumber\\
 &\equiv& k_1+k_2+k_3+k_4.
\end{eqnarray}
First, we consider $k_1$. Then by Lemma 3.9, it holds uniformly for all $y_1\in [0,a]$ that
 \begin{equation}J_1\sim n\overline{F}_1 (u_1),\ {\rm as}\  u_1\to\infty \end{equation}
and it holds uniformly for all $y_2\in [0,b]$ that
 \begin{equation}
 J_2\sim n\overline{F}_2 (u_2),\ {\rm as}\  u_2\to\infty.
 \end{equation}
Using Lemma 3.1 and the  dominated convergence theorem, we get
\begin{equation}
k_1\sim n^2 \overline{F}_1 (u_1)\overline{F}_2 (u_2) \int_0^{a}\int_0^{b}P(\overline{B_1}(T)\in dy_1, \overline{B_2}(T)\in dy_2), \ {\rm as}\  (u_1, u_2)\to (\infty,\infty).\nonumber
\end{equation}
Thus
\begin{equation}
\lim_{(a,b)\to (\infty,\infty)}\lim_{(u_1, u_2)\to (\infty,\infty)}\frac{k_1}{n^2 \overline{F}_1 (u_1)\overline{F}_2 (u_2)}=1.
\end{equation}
Now consider $k_2$. Using (4.11), Lemma 3.1 and the  dominated convergence theorem,
\begin{eqnarray}
k_2 &\sim & n \overline{F}_1 (u_1)\int_0^{a}\int_b^{\infty}P(\overline{B_1}(T)\in dy_1, \overline{B_2}(T)\in dy_2)J_2\nonumber\\
&\le & n \overline{F}_1 (u_1) P\left(\sum_{i=1}^{n}X_{2i}-C_2(T)-\sigma_2 b>u_2\right)\int_0^{a}\int_b^{\infty}P(\overline{B_1}(T)\in dy_1, \overline{B_2}(T)\in dy_2) \nonumber\\
&\sim & n^2 \overline{F}_1 (u_1)\overline{F}_2 (u_2) \int_0^{a}\int_{b}^{\infty}P(\overline{B_1}(T)\in dy_1, \overline{B_2}(T)\in dy_2)
 ,\ {\rm as}\  (u_1, u_2)\to (\infty,\infty).\nonumber
\end{eqnarray}
Thus
\begin{equation}
\lim_{(a,b)\to (\infty,\infty)}\lim_{(u_1, u_2)\to (\infty,\infty)}\frac{k_2}{n^2 \overline{F}_1 (u_1)\overline{F}_2 (u_2)}=0.
\end{equation}
Likewise,
 \begin{equation}
\lim_{(a,b)\to (\infty,\infty)}\lim_{(u_1, u_2)\to (\infty,\infty)}\frac{k_3}{n^2 \overline{F}_1 (u_1)\overline{F}_2 (u_2)}=0.
\end{equation}
Finally we deal with $k_4$.
\begin{eqnarray}
k_4 & \le & P\left(\sum_{i=1}^{n}X_{1i}-C_1(T)-\sigma_1 a>u_1\right)P\left(\sum_{i=1}^{n}X_{2i}-C_2(T)-\sigma_2 b>u_2\right)\nonumber\\
&&\times\int_{a}^{\infty}\int_{b}^{\infty}P(\overline{B_1}(T)\in dy_1, \overline{B_2}(T)\in dy_2)\nonumber\\
& \sim & n^2 \overline{F}_1 (u_1)\overline{F}_2 (u_2)\int_{a}^{\infty}\int_{b}^{\infty}P(\overline{B_1}(T)\in dy_1, \overline{B_2}(T)\in dy_2),\ {\rm as}\  (u_1, u_2)\to (\infty,\infty),\nonumber
\end{eqnarray}
from which we get
 \begin{equation}
\lim_{(a,b)\to (\infty,\infty)}\lim_{(u_1, u_2)\to (\infty,\infty)}\frac{k_4}{n^2 \overline{F}_1 (u_1)\overline{F}_2 (u_2)}=0.
\end{equation}
From (4.9) and (4.13)-(4.16) we obtain
\begin{equation}
 \lim_{(u_1, u_2)\to (\infty,\infty)}\frac{I_1}{n^2 \overline{F}_1 (u_1)\overline{F}_2 (u_2)}=1.
\end{equation}
Now, it follows from (4.8), (4.17) and the dominated convergence theorem, we find that
\begin{equation}
 \lim_{(u_1, u_2)\to (\infty,\infty)}\frac{\psi_{\max}(\vec{u};T)}{\lambda T(1+\lambda T)\overline{F_1}(u_1)\overline{F_2}(u_2)}\ge 1,\nonumber
\end{equation}
from which and (4.7), we obtain (4.1).

Note that
$$\psi_{\rm{and}}(\vec{u};T)\le P\left(\sum_{i=1}^{N(T)}X_{1i}-\sigma_1 \underline{B}_1(T)>u_1, \sum_{i=1}^{N(T)}X_{2i}-\sigma_2 \underline{B}_2(T)>u_2\right),$$
from which, (4.4) and (4.7) we have
$$\lim_{(u_1, u_2)\to (\infty,\infty)}\frac{\psi_{\rm{and}}(\vec{u};T)}{\overline{F_1}(u_1)+\overline{F_2}(u_2)}\le\lim_{(u_1, u_2)\to (\infty,\infty)}\frac{\lambda T(1+\lambda T)\overline{F_1}(u_1)\overline{F_2}(u_2)}{\overline{F_1}(u_1)+\overline{F_2}(u_2)}=0.$$
\textcolor{blue} {Thus, $\psi_{and}(\vec{u};T)\sim 0$, as $(u_1, u_2) \rightarrow (\infty, \infty).$ By (1.6), we have that
$$\psi_{min}(\vec{u};T)\sim P(T_1\leq T\mid U_1(0)=u_1)+P(T_2\leq T\mid U_2(0)=u_1)=\psi_1(u_1;T)+\psi_2(u_2;T).$$
From Lemma 3.3, we can obtain (4.2).}

Next, we prove relation (4.3). Using Theorem 7.2 in Ikeda and Watanabe (1981, P. 85) (see also Yin and Wen (2013)) one has, for all $t\ge 0$,
$$\sqrt{\sigma_1^2+\sigma_2^2+2\rho\sigma_1\sigma_2}W(t)\stackrel{d}{=}\sigma_1 B_1(t)+\sigma_2 B_2(t),$$
where `$\stackrel{d}{=}$' denotes equality in distribution, $W$ is a standard Brownian motion independent of $\{X_{1k}, k=1,2,\cdots\}$, $\{X_{2k}, k=1,2,\cdots\}$, $\{C_1(t), t\ge 0\}$, $\{C_2(t), t\ge 0\}$ and  $\{N(t), t\ge 0\}$.
Thus, for all $t\ge 0$,  $U_1(t)+U_2(t)$ can be written as
$$U_1(t)+U_2(t)\stackrel{d}{=}u_1+u_2+C_1(t)+C_2(t)-\sum_{i=1}^{N(t)}(X_{1i}+X_{2i})
+\sqrt{\sigma_1^2+\sigma_2^2+2\rho\sigma_1\sigma_2}W(t).$$
Applying Lemma 3.3 to this model, we get that
if $F_1*F_2\in\mathcal{S}$, then
$$\psi_{\rm{sum}}(\vec{u}; T)\sim  \lambda T\overline{F_1*F_2}(u_1+u_2)
\sim  \lambda T(\overline{F_1}(u_1+u_2)+\overline{F_2}(u_1+u_2)), \ u_1+u_2\to\infty,$$
where, in the last step, we have used the fact in  Embrechts and Goldie (1980) (see also Geluk and Tang (2009)) which states that
$$F_1*F_2\in\mathcal{S} \ {\rm if\ and\ only\ if}\ P(X_1+X_2>x)\sim
\overline{F_1}(x)+\overline{F_2}(x).$$
This ends the proof of Theorem 4.1.
\begin{remark} Letting  $\{C_i(t)=c_i t, i=1,2$ and $\rho=0$ in Theorem 4.1, we get Theorem 4.1 in Li, Liu and Tang (2007).
\end{remark}

\begin{theorem} Consider the insurance risk model introduced in Section 1. Assume that $N_1(t)=N_2(t)=N(t)$, $\rho\in (-1,0]$, $r>0$ and $\{X_{1k}, k=1,2,\cdots\}$, $\{X_{2k}, k=1,2,\cdots\}$, $\{C_1(t), t\ge 0\}$, $\{C_2(t), t\ge 0\}$,  $\{N(t), t\ge 0\}$, $\{(B_1(t), B_2(t)), t\ge 0\}$
are mutually independent. \\
(a)\ If  $F_1, F_2\in \mathcal{S}$, then, for each fixed time $T\ge 0$,
\begin{equation}
\psi_{\max}(\vec{u};T)\sim\frac{\lambda(\lambda+\frac{1}{T})}{r^2}\int_{u_1}^{u_1 e^{rT}}\frac{\overline{F_1}(y)}{y}dy
\int_{u_2}^{u_2 e^{rT}}\frac{\overline{F_2}(y)}{y}dy,\ \  \text{\rm as}\ \ (u_1,u_2)\to (\infty,\infty),
 \end{equation}
 \begin{equation}
 \psi_{\min}(\vec{u};T)\sim \frac{\lambda}{r}\left(\int_{u_1}^{u_1 e^{rT}}\frac{\overline{F_1}(y)}{y}dy+
\int_{u_2}^{u_2 e^{rT}}\frac{\overline{F_2}(y)}{y}dy\right),\ \  \text{\rm as}\ \ (u_1,u_2)\to (\infty,\infty).
 \end{equation}
 (b)\ If  $F_1*F_2\in \mathcal{S}$, then, for each fixed time $T\ge 0$,
 \begin{equation}
 \psi_{\rm{sum}}(\vec{u};T)\sim \lambda T\int_0^1\overline{F_1*F_2}(e^{rTz}(u_1+u_2))dz,\ \  \text{\rm as}\ \ u_1+u_2\to \infty.
 \end{equation}
 In particular, if there are two positive constants $l_1$ and $l_2$ such that
 $\overline{F}_i (x)\sim l_i \overline{F} (x)$, $i=1,2$. Then
 \begin{equation}
 \psi_{\rm{sum}}(\vec{u};T)\sim \lambda T\left(\int_0^1\overline{F_1}(e^{rTz}(u_1+u_2))+\int_0^1\overline{F_2}(e^{rTz}(u_1+u_2))\right),\ \  \text{\rm as}\ \ u_1+u_2\to \infty.
 \end{equation}
\end{theorem}
{\bf Proof.}\ We write $\psi_{\max}(\vec{u};T)$ as
$$\psi_{\max}(\vec{u};T)=P(e^{-rt} U_i(t)<0, i=1,2 \ {\rm for}\ {\rm some}\  0<t\le T|\vec{U}(0)=\vec{u}).$$
For $t\in [0,T]$ and each $i=1$ or $2$, we have
 \begin{eqnarray}
 u_i-\int_0^t e^{-rs}dS_i(s)+\sigma_i \int_0^t e^{-rs}dB_i(s)&\le& e^{-rt} U_i(t)
 \le u_i+\int_0^T  e^{-rs}dC_i(s)\nonumber\\
 &&-\int_0^t e^{-rs}dS_i(s)
+\sigma_i \int_0^t e^{-rs}dB_i(s).\nonumber
\end{eqnarray}
It follows that $\psi_{\max}(\vec{u};T)$ satisfies
 \begin{eqnarray}
 \psi_{\max}(\vec{u};T)&\le& P\left(\sum_{i=1}^{N(T)}\vec{X}_i e^{-r\tau_i}-{\sigma_1  \underline{\Delta}_1(T)\choose
  \sigma_2  \underline{\Delta}_2(T)}>\vec{u}\right))\nonumber\\
 &\le&\sum_{n=0}^{\infty}P(N(T)=n)P\left(\sum_{i=1}^{n}\vec{X}_i e^{-r\tau_i}-{\sigma_1 \underline{\Delta}_1(T)\choose
  \sigma_2 \underline{\Delta}_2(T)}>\vec{u}\big|N(t)=n\right)\nonumber\\
  &\le&\sum_{n=0}^{\infty}P(N(T)=n)\int_0^{\infty}\int_0^{\infty}P\left(\sum_{i=1}^{n}\vec{X}_i e^{-rT U_i}\in d\vec{z}\right)\nonumber\\
  &&\times P\left(\vec{z}-{\sigma_1 \underline{\Delta}_1(T)\choose
  \sigma_2 \underline{\Delta}_2(T)}>\vec{u}\right).
\end{eqnarray}
where we have used Lemma 3.7 in the last steps.
Since $\rho\in (-1,0]$, using Lemma 3.6 one has
\begin{eqnarray}
&&P\left(\vec{z}-{\sigma_1 \underline{\Delta}_1(T)\choose
\sigma_2 \underline{\Delta}_2(T)}>\vec{u}\right)\nonumber\\
&&\le P(z_1-\sigma_1 \underline{\Delta}_1(T)>u_1)P(z_2-\sigma_2\underline{\Delta}_2(T)>u_2).
\end{eqnarray}
Using independence of $\{X_{1k}, k=1,2,\cdots\}$ and $\{X_{2k}, k=1,2,\cdots\}$, we have
\begin{eqnarray}
P\left(\sum_{i=1}^{n}\vec{X}_i e^{-rT U_i}\in d\vec{z}\right)&=&\int_0^1\cdots\int_0^1
P\left(\sum_{i=1}^{n}X_{1i} e^{-rT v_i}\in d z_1\right)P\left(\sum_{i=1}^{n}X_{2i} e^{-rT v_i}\in d z_2\right)\nonumber\\
&&\times\prod_{j=1}^n P(U_j\in d v_j).
\end{eqnarray}
Substituting (4.23) and (4.24) into (4.22), and using the following
$$ P\left( \sum_{i=1}^{n}X_{1i} e^{-rT v_i}-\sigma_1 \underline{\Delta}_1(T)>u_1\right)\sim
 P\left( \sum_{i=1}^{n}X_{1i} e^{-rT v_i}>u_1\right),  u_1\to\infty,
$$
and
$$ P\left( \sum_{i=1}^{n}X_{2i} e^{-rT v_i}-\sigma_2 \underline{\Delta}_2(T)>u_2\right)\sim
 P\left( \sum_{i=1}^{n}X_{2i} e^{-rT v_i}>u_2\right),  u_2\to\infty,
$$
uniformly for $(v_1,\cdots,v_n)\in [0,1]^n$,
we obtain
 \begin{eqnarray}
\psi_{\max}(\vec{u};T)&\lesssim& \sum_{n=0}^{\infty}P(N(T)=n) P\left(\sum_{i=1}^{n}X_{1i} e^{-rTU_i}>u_1, \sum_{i=1}^{n}X_{2i} e^{-rTU_i}>u_2\right)\nonumber\\
&\equiv& \sum_{n=0}^{\infty}P(N(T)=n) k_5.
\end{eqnarray}
We apply Proposition 5.1 of Tang and Tsitsiashvili (2003), which says that, for i.i.d. subexponential random variables $\{X_k\}$ and for arbitrarily $a$ and $b$, $0<a\le b<\infty$, the relation
$$  P\left(\sum_{i=1}^{n}c_iX_i >x\right)\sim \sum_{i=1}^n P( c_iX_i  >x)$$ holds uniformly for $(c_1,\cdots,c_n)\in [a,b]\times\cdots\times[a,b]$. Hence, by conditioning on $(U_1, \cdots, U_n),$ we find that
where
\begin{eqnarray}
 k_5\sim n^2  P\left(X_{11} e^{-rTU_1}>u_1\right) P\left(X_{21} e^{-rTU_1}>u_2\right).
\end{eqnarray}
Substituting (4.26) into (4.25), and using  the  dominated convergence theorem, we get
\begin{eqnarray}
\limsup_{(u_i,u_2)\to (\infty,\infty)}\frac{\psi_{\max}(\vec{u};T)}{\frac{\lambda(\lambda+\frac{1}{T})}{r^2}\int_{u_1}^{u_1 e^{rT}}\frac{\overline{F_1}(y)}{y}dy
\int_{u_2}^{u_2 e^{rT}}\frac{\overline{F_2}(y)}{y}dy}\le 1.
\end{eqnarray}
Next,  we establish the asymptotic lower bound for $\psi_{\max}(\vec{u};T)$.
Clearly,
 \begin{eqnarray}
 \psi_{\max}(\vec{u};T)&\ge& P\left(\sum_{i=1}^{N(T)}\vec{X}_i e^{-r\tau_i}-\int_0^T e^{-rs}d\vec{C}(s)-{\sigma_1 \overline{\Delta}_1(T)\choose
  \sigma_2 \overline{\Delta}_2(T)}>\vec{u}\right) \nonumber\\
 &=&\sum_{n=0}^{\infty}P(N(T)=n)P\left(\sum_{i=1}^{n}\vec{X}_i e^{-rTU_i}- {\sigma_1 \overline{\Delta}_1(T)\choose
  \sigma_2 \overline{\Delta}_2(T)}-\int_0^T e^{-rs}d\vec{C}(s)>\vec{u}\right)\nonumber\\
  &\equiv& \sum_{n=0}^{\infty}P(N(T)=n) I_2,
\end{eqnarray}
where, for some positive constants $c$ and $d$,
$$I_2=\left(\int_0^{c}\int_0^{d}+\int_0^{c}\int_d^{\infty}+\int_c^{\infty}\int_0^{d}
+\int_c^{\infty}\int_d^{\infty}\right)P(\overline{\Delta_1}(T)\in dy_1, \overline{\Delta_2}(T)\in dy_2)J_3 J_4.$$
Here,
$$J_3=P\left(\sum_{i=1}^{n}X_{1i}e^{-rTU_i}-\int_0^T e^{-rs} dC_1(s)-\sigma_1 y_1>u_1\right),$$
and
$$J_4=P\left(\sum_{i=1}^{n}X_{2i}e^{-rTU_i}-\int_0^T e^{-rs} dC_2(s) -\sigma_2 y_2>u_2\right).$$
By Lemma 3.8 we know that
$\sum_{i=1}^{n}X_{ji}e^{-rTU_i}\in \mathcal{S}, j=1,2$, since all $X_{ji}\in \mathcal{S}$.
Then invoke Lemma 3.9, we get
$$J_3\sim n P(X_{11}e^{-rTU_1}>u_1),\ {\rm as}\ u_1\to\infty, \ J_4\sim n P(X_{21}e^{-rTU_1}>u_2)\ {\rm as}\ u_2\to\infty$$
uniformly for all $y_1\in [0,c]$ and $ y_2\in [0,d]$, respectively.
Now, using the same argument as leads to (4.17), we have
\begin{equation}
 \lim_{(u_1, u_2)\to (\infty,\infty)}\frac{I_2}{n^2  P(X_{11}e^{-rTU_1}>u_1) P(X_{21}e^{-rTU_1}>u_2)}=1.
\end{equation}
Now, it follows from (4.28), (4.29), Lemma 3.1 and the dominated convergence theorem that
$$
\lim_{(u_i,u_2)\to (\infty,\infty)}\frac{\psi_{\max}(\vec{u};T)}{\lambda T(1+\lambda T) P(X_{11}e^{-rTU_1}>u_1) P(X_{21}e^{-rTU_1}>u_2)}\ge 1,
 $$
or, equivalently,
  $$
\lim_{(u_i,u_2)\to (\infty,\infty)}\frac{\psi_{\max}(\vec{u};T)}{\frac{\lambda(\lambda+\frac{1}{T})}{r^2}\int_{u_1}^{u_1 e^{rT}}\frac{\overline{F_1}(y)}{y}dy
\int_{u_2}^{u_2 e^{rT}}\frac{\overline{F_2}(y)}{y}dy}\ge 1,
 $$
 from which and (4.27), we obtain (4.18).

 The relation (4.19) follows from (1.6) and Lemma 3.4 since, as above,
 \begin{eqnarray}
 &\lim\limits_{(u_1, u_2)\to (\infty,\infty)}&\frac{\psi_{\rm{and}}(\vec{u};T)}
{\int_{u_1}^{u_1 e^{rT}}\frac{\overline{F_1}(y)}{y}dy+
\int_{u_2}^{u_2 e^{rT}}\frac{\overline{F_2}(y)}{y}dy}\nonumber\\
&&\le\lim_{(u_1, u_2)\to (\infty,\infty)}\frac {\frac{\lambda(\lambda+\frac{1}{T})}{r^2}\int_{u_1}^{u_1 e^{rT}}\frac{\overline{F_1}(y)}{y}dy
\int_{u_2}^{u_2 e^{rT}}\frac{\overline{F_2}(y)}{y}dy}
{\int_{u_1}^{u_1 e^{rT}}\frac{\overline{F_1}(y)}{y}dy+
\int_{u_2}^{u_2 e^{rT}}\frac{\overline{F_2}(y)}{y}dy}=0. \nonumber
\end{eqnarray}
\textcolor{blue} { By (1.6), we have that
$$\psi_{min}(\vec{u};T)\sim \psi_1(u_1;T)+\psi_2(u_2;T),\ \ \text{\rm as} \ \ (u_1, u_2) \rightarrow (\infty, \infty).$$
From Lemma 3.4, we have that
$$\psi_i(u_i; T)\sim
\frac{\lambda}{r}\int_{u_i}^{u_i e^{rT}}\frac{\overline{F_i}(y)}{y}dy,\  u_i\to\infty,\ i=1,2.$$ }
\textcolor{blue}{ Then,
$$ \psi_{min}(\vec{u};T)\sim \frac{\lambda}{r}\left(\int_{u_1}^{u_1 e^{rT}}\frac{\overline{F_1}(y)}{y}dy+
\int_{u_2}^{u_2 e^{rT}}\frac{\overline{F_2}(y)}{y}dy\right),\ \  \text{\rm as}\ \ (u_1,u_2)\to (\infty,\infty).$$
Thus, we complete the proof of (4.19).}

Next, we prove relation (4.20). Similarly, we have, for all $t\ge 0$,
\begin{eqnarray}
U_1(t)+U_2(t)&\stackrel{d}{=}&(u_1+u_2)e^{rt}+\int_0^t e^{r(t-s)}d(C_1(s)+C_2(s))\nonumber\\
&&-\int_0^t e^{r(t-s)}d\sum_{i=1}^{N(s)}(X_{1i}+X_{2i})\nonumber\\
&&+\sqrt{\sigma_1^2+\sigma_2^2+2\rho\sigma_1\sigma_2}\int_0^t e^{r(t-s)}dW(s),
\end{eqnarray}
where $\{W(t), t\ge 0\}$ is a standard Brownian motion independent of $\{X_{1k}, k=1,2,\cdots\}$, $\{X_{2k}, k=1,2,\cdots\}$, $\{C_1(t), t\ge 0\}$, $\{C_2(t), t\ge 0\}$ and  $\{N(t), t\ge 0\}$.

\textcolor{blue}{From Lemma 3.4 we have that
$$\psi_{sum}(\vec{u}; T)\sim
\frac{\lambda}{r}\int_{u_1+u_2}^{(u_1+u_2) e^{rT}}\frac{\overline{F_1\ast F_2}(y)}{y}dy,\  u_1+u_2\to\infty.$$
Let $y=(u_1+u_2) e^{rTz}$, then $dy=rT(u_1+u_2) e^{rTz}dz$. Therefore,}
\textcolor{blue}{\begin{align}
\nonumber \psi_{sum}(\vec{u}; T)&\sim
\frac{\lambda}{r}\int_{0}^{1}\frac{\overline{F_1\ast F_2}((u_1+u_2) e^{rTz})}{(u_1+u_2) e^{rTz}}rT(u_1+u_2) e^{rTz}dz\\
\nonumber&=T\lambda\int_{0}^{1}\overline{F_1\ast F_2}((u_1+u_2) e^{rTz})dz,  \text{\rm as}\ \ (u_1,u_2)\to (\infty,\infty).
\end{align}
This completes the proof of (4.20).}
The result (4.21) follows from (4.20) and Lemma 3.1 in Hao and Tang (2008).
This ends the proof of Theorem 4.2.
\begin{remark} Letting  $\{C_i(t)=c_i t, i=1,2$,  $\rho=0$  $\sigma_1=0, \sigma_2=0$  in Theorem 4.2, we get the result in Liu, Wang and Long (2007).
\end{remark}

\begin{theorem} Consider the insurance risk model introduced in Section 1. Assume that  $\rho\in (-1,0]$, $r=0$ and $\{X_{1k}, k=1,2,\cdots\}$, $\{X_{2k}, k=1,2,\cdots\}$, $\{C_1(t), t\ge 0\}$, $\{C_2(t), t\ge 0\}$,  $\{N_i(t), t\ge 0\}$, $i=1,2$, $\{(B_1(t), B_2(t)), t\ge 0\}$
are mutually independent.\\
(a)\  If  $F_1, F_2\in \mathcal{S}$, then, for each fixed time $T\ge 0$,
\begin{equation}
\psi_{\max}(\vec{u};T)\sim \lambda_1 \lambda_2 T^2\overline{F_1}(u_1)\overline{F_2}(u_2),\ \  \text{\rm as}\ \ (u_1,u_2)\to (\infty,\infty),
\end{equation}
\begin{equation}
\psi_{\min}(\vec{u};T)\sim  T\left(\lambda_1\overline{F_1}(u_1)+\lambda_2\overline{F_2}(u_2)\right),\ \  \text{\rm as}\ \ (u_1,u_2)\to (\infty,\infty).
\end{equation}
 (b)\ If  $F_{\xi X_{11}+(1-\xi)X_{21}}\in \mathcal{S}$, where $\xi$ is a random variable  independent of $\{X_{1k}, k=1,2,\cdots\}$ and $\{X_{2k}, k=1,2,\cdots\}$, and $P(\xi=1)=1-P(\xi=0)=\frac{\lambda_1}{\lambda_1+\lambda_2}$,  then, for each fixed time $T\ge 0$,
 \begin{equation}
 \psi_{\rm{sum}}(\vec{u};T)\sim T\left( \lambda_1\overline{F_1}(u_1+u_2)+ \lambda_2 \overline{F_2}(u_1+u_2)\right),\ \  \text{\rm as}\ \ u_1+u_2\to \infty.
 \end{equation}
\end{theorem}

{\bf Proof.}\  The proof is similar to that of Theorem 4.1, so we only give the main steps. First we establish the asymptotic upper bound for $\psi_{\max}(\vec{u};T)$.
Clearly,
 \begin{eqnarray}
 \psi_{\max}(\vec{u};T)&\le& P\left({\sum_{i=1}^{N_1(T)}X_{1i}\choose \sum_{i=1}^{N_2(T)}X_{2i}}-{\sigma_1 \underline{B}_1(T)\choose \sigma_2 \underline{B}_2(T)}> {u_1\choose u_2}\right)\nonumber\\
  &=&\int_0^{\infty}\int_0^{\infty}P\left(\sum_{i=1}^{N_1(T)}X_{1i}\in d z_1\right) P\left(\sum_{i=1}^{N_2(T)}X_{2i}\in d z_2\right)\nonumber\\
  &&\times P\left({z_1\choose z_2}-{\sigma_1 \underline{B}_1(T)\choose \sigma_2 \underline{B}_2(T)}>{u_1\choose u_2}\right).
\end{eqnarray}
Since $\rho\in (-1,0]$, using (3.8) one has
\begin{eqnarray}
&P&\left({z_1\choose z_2}-{\sigma_1 \underline{B}_1(T)\choose \sigma_2 \underline{B}_2(T)}>{u_1\choose u_2}\right)\nonumber\\
&&\le P(z_1-\sigma_1 \underline{B}_1(T)>u_1)P(z_2-\sigma_2\underline{B}_2(T)>u_2).
\end{eqnarray}
Substituting (4.35)  into (4.34) we get
 \begin{eqnarray}
\psi_{\max}(\vec{u};T)&\le&   P\left(\sum_{i=1}^{N_1(T)}X_{1i}-\sigma_1 \underline{B}_1(T)>u_1\right)
  P\left(\sum_{i=1}^{N_2(T)}X_{2i}-\sigma_2 \underline{B}_2(T)>u_1\right)\nonumber\\
 &\sim& \lambda_1\lambda_2 T^2\overline{F_1}(u_1)\overline{F_2}(u_2),\ \  \text{\rm as}\ \ (u_1,u_2)\to (\infty,\infty),
 \end{eqnarray}
where in the last step we have used Lemma 3.3.

Next,  we establish the asymptotic lower bound for $\psi_{\max}(\vec{u};T)$.
Clearly,
 \begin{eqnarray}
 \psi_{\max}(\vec{u};T)&\ge& P\left({\sum_{i=1}^{N_1(T)}X_{1i}\choose \sum_{i=1}^{N_2(T)}X_{2i}}-{C_1(T)\choose C_2(T)}- {\sigma_1 \overline{B}_1(T)\choose \sigma_2 \overline{B}_2(T)}>{u_1\choose u_2}\right)\nonumber\\
  &=& \sum_{n=0}^{\infty}P(N_1(T)=n) \sum_{m=0}^{\infty}P(N_1(T)=m)I_3,
\end{eqnarray}
where
$$I_3=P\left({\sum_{i=1}^{n}X_{1i}\choose \sum_{i=1}^{m}X_{2i}}-{C_1(T)\choose C_2(T)}- {\sigma_1 \overline{B}_1(T)\choose \sigma_2 \overline{B}_2(T)}>{u_1\choose u_2}\right).$$
Using the same arguments as proving (4.17), we get
\begin{equation}
 \lim_{(u_1, u_2)\to (\infty,\infty)}\frac{I_3}{nm \overline{F}_1 (u_1)\overline{F}_2 (u_2)}=1,\nonumber
\end{equation}
which, together with (4.37), we have
\begin{equation}
 \lim_{(u_1, u_2)\to (\infty,\infty)}\frac{\psi_{\max}(\vec{u};T)}{\lambda_1 \lambda_2 T^2\overline{F_1}(u_1)\overline{F_2}(u_2)}\ge 1.\nonumber
\end{equation}
The proof of (4.32) is straightforward and thus we omit it. Next, we prove (4.33). Using the properties of two independent compound Poisson processes and two independent  Brownian motions, we have, for all $t\ge 0$,
\begin{eqnarray}
U_1(t)+U_2(t)&\stackrel{d}{=}&u_1+u_2+C_1(t)+C_2(t)-\sum_{i=1}^{N_0(t)}(\xi X_{1i}+(1-\xi)X_{2i})\nonumber\\
&&+\sqrt{\sigma_1^2+\sigma_2^2+2\rho\sigma_1\sigma_2}W(t),\nonumber
\end{eqnarray}
where $\{W(t), t\ge 0\}$ is a standard Brownian motion, $\{N_0(t), t\ge 0\}$ is a Poisson process with intensity $\lambda_1+\lambda_2$,  $\xi$ is a Bernoulli  random variable with $P(\xi=1)=1-P(\xi=0)=\frac{\lambda_1}{\lambda_1+\lambda_2}$. Moreover, $\xi $, $\{W(t), t\ge 0\}$, $\{N_0(t), t\ge 0\}$, $\{X_{1k}, k=1,2,\cdots\}$, $\{X_{2k}, k=1,2,\cdots\}$, $\{C_1(t), t\ge 0\}$, $\{C_2(t), t\ge 0\}$ and  $\{N(t), t\ge 0\}$ are independent.
Applying Lemma 3.3 to this model, we get
$$\psi_{\rm{sum}}(\vec{u}; T)\sim  (\lambda_1+\lambda_2) T\overline{F}_{\xi X_{11}+(1-\xi) X_{21}}(u_1+u_2), \ u_1+u_2\to\infty,$$ and  the result (4.33) follows since (c.f. Kaas et al. 2008)
$$P(\xi X_{11}+(1-\xi) X_{21}>u_1+u_2)=\frac{\lambda_1}{\lambda_1+\lambda_2}\overline{F}_1(u_1+u_2)+ \frac{\lambda_2}{\lambda_1+\lambda_2}\overline{F}_2(u_1+u_2).$$
This ends the proof of Theorem 4.3.

\begin{theorem} Consider the insurance risk model introduced in Section 1. Assume that  $\rho\in (-1,0]$, $r>0$ and $\{X_{1k}, k=1,2,\cdots\}$, $\{X_{2k}, k=1,2,\cdots\}$, $\{C_1(t), t\ge 0\}$, $\{C_2(t), t\ge 0\}$,  $\{N_i(t), t\ge 0\}$, $i=1,2$, $\{(B_1(t), B_2(t)), t\ge 0\}$
are mutually independent. \\
(a)\ If  $F_1, F_2\in \mathcal{S}$, then, for each fixed time $T\ge 0$,
\begin{equation}
\psi_{\max}(\vec{u};T)\sim\frac{\lambda_1\lambda_2}{r^2}\int_{u_1}^{u_1 e^{rT}}\frac{\overline{F_1}(y)}{y}dy
\int_{u_2}^{u_2 e^{rT}}\frac{\overline{F_2}(y)}{y}dy,\ \  \text{\rm as}\ \ (u_1,u_2)\to (\infty,\infty),
 \end{equation}
 \begin{equation}
 \psi_{\min}(\vec{u};T)\sim \frac{1}{r}\left(\lambda_1\int_{u_1}^{u_1 e^{rT}}\frac{\overline{F_1}(y)}{y}dy+
\lambda_2\int_{u_2}^{u_2 e^{rT}}\frac{\overline{F_2}(y)}{y}dy\right),\ \  \text{\rm as}\ \ (u_1,u_2)\to (\infty,\infty).
 \end{equation}
 (b)\ If  $F_{\xi X_{11}+(1-\xi)X_{21}}\in \mathcal{S}$, where $\xi$ is defined as in Theorem 4.3, then, for each fixed time $T\ge 0$,
 \begin{equation}
 \psi_{\rm{sum}}(\vec{u};T)\sim  \frac{1}{r}\left(\lambda_1\int_{u_1+u_2}^{(u_1+u_2) e^{rT}}\frac{\overline{F_1}(y)}{y}dy+
\lambda_2\int_{u_1+u_2}^{(u_1+u_2) e^{rT}}\frac{\overline{F_2}(y)}{y}dy\right),\ \  \text{\rm as}\ \ u_1+u_2\to \infty.
 \end{equation}
\end{theorem}

{\bf Proof.}\  As in the proof of Theorem 4.2, we have
 \begin{eqnarray}
 \psi_{\max}(\vec{u};T)&\le& P\left({\sum_{i=1}^{N_1(T)}X_{1i}e^{-r\tau_i}\choose \sum_{i=1}^{N_2(T)}X_{2i}e^{-r\tau_i}}-{\sigma_1  \underline{\Delta}_1(T)\choose
  \sigma_2  \underline{\Delta}_2(T)}>\vec{u}\right)\nonumber\\
 &\le&\sum_{n=0}^{\infty}P(N_1(T)=n)\sum_{m=0}^{\infty}P(N_2(T)=m)\nonumber\\
  &&\times P\left( {\sum_{i=1}^{n}X_{1i}e^{-r\tau_i}\choose \sum_{i=1}^{m}X_{2i}e^{-r\tau_i}}-{\sigma_1 \underline{\Delta}_1(T)\choose
  \sigma_2 \underline{\Delta}_2(T)}>\vec{u}\right)\nonumber\\
  &\lesssim&\sum_{n=0}^{\infty}\sum_{m=0}^{\infty}  nm P(N(T)=n)P(N_2(T)=m)
   P\left(X_{11} e^{-rTU_1}>u_1\right) P\left(X_{21} e^{-rTU_1}>u_2\right)\nonumber\\
   &=&\lambda_1\lambda_2 T^2  P\left(X_{11} e^{-rTU_1}>u_1\right) P\left(X_{21} e^{-rTU_1}>u_2\right).\nonumber
\end{eqnarray}
It follows that
\begin{eqnarray}
\limsup_{(u_i,u_2)\to (\infty,\infty)}\frac{\psi_{\max}(\vec{u};T)}{\frac{\lambda_1\lambda_2}{r^2}\int_{u_1}^{u_1 e^{rT}}\frac{\overline{F_1}(y)}{y}dy
\int_{u_2}^{u_2 e^{rT}}\frac{\overline{F_2}(y)}{y}dy}\le 1.\nonumber
\end{eqnarray}
The asymptotic lower bound for $\psi_{\max}(\vec{u};T)$ can be established similarly.

The relation (4.39) follows from (1.6), Lemma 3.4 and the fact
\begin{eqnarray}
&\lim_{(u_1, u_2)\to (\infty,\infty)}&\frac{\psi_{\rm{and}}(\vec{u};T)}
{\lambda_1\int_{u_1}^{u_1 e^{rT}}\frac{\overline{F_1}(y)}{y}dy+\lambda_2
\int_{u_2}^{u_2 e^{rT}}\frac{\overline{F_2}(y)}{y}dy}\nonumber\\
&&\le \lim_{(u_1, u_2)\to (\infty,\infty)}\frac{\frac{\lambda_1\lambda_2}{r^2}\int_{u_1}^{u_1 e^{rT}}\frac{\overline{F_1}(y)}{y}dy
\int_{u_2}^{u_2 e^{rT}}\frac{\overline{F_2}(y)}{y}dy}
{\lambda_1\int_{u_1}^{u_1 e^{rT}}\frac{\overline{F_1}(y)}{y}dy+\lambda_2
\int_{u_2}^{u_2 e^{rT}}\frac{\overline{F_2}(y)}{y}dy}=0.\nonumber
\end{eqnarray}
Finally, we prove (4.40). Using the same arguments as above, we have
\begin{eqnarray}
U_1(t)+U_2(t)&\stackrel{d}{=}&(u_1+u_2)e^{rt}+\int_0^t e^{r(t-s)}d(C_1(s)+C_2(s))\nonumber\\
&&-\int_0^t e^{r(t-s)}d \sum_{i=1}^{N_0(t)}(\xi X_{1i}+(1-\xi)X_{2i})\nonumber\\
&&+\sqrt{\sigma_1^2+\sigma_2^2+2\rho\sigma_1\sigma_2}\int_0^t e^{r(t-s)}dW(s),\ t\ge 0,
\end{eqnarray}
where $\xi $, $\{W(t), t\ge 0\}$, $\{N_0(t), t\ge 0\}$ are the same as in the proof of Theorem 4.3.
It follows from Lemma 3.4 that
$$\psi_{\rm{sum}}(\vec{u}; T)\sim  \frac{\lambda_1+\lambda_2}{r}\int_{u_1+u_2}^{(u_1+u_2) e^{rT}}
\frac{\overline{F}_{\xi X_{11}+(1-\xi) X_{21}}(y)}{y}dy, \ u_1+u_2\to\infty,$$
and  the result (4.40) follows since
$$\overline{F}_{\xi X_{11}+(1-\xi) X_{21}}(y) =\frac{\lambda_1}{\lambda_1+\lambda_2}\overline{F}_1(y)+ \frac{\lambda_2}{\lambda_1+\lambda_2}\overline{F}_2(y).$$
This completes the proof of Theorem 4.4.

\textcolor{blue} {\section{ Conclusion}}

\textcolor{blue} {In this paper, we investigated the bidimensional risk model that describe the surplus process of an insurer. We give new results for the different types of finite-time ruin probabilities under the circumstance of that the Brownian motions are correlated with constant correlation coefficient.  We remark that the extension to multidimensional models is more complicated. However, the multidimensional model can better describe different insurance businesses. In addition, we can also consider the relationship between different businesses in the future research, which I believe that can be a more interesting problem.}
\setcounter{equation}{0}

\noindent{\bf Acknowledgements.}

The research was supported by the National Natural Science Foundation of China (No. 12071251); the Youth Innovation Team of Shandong Universities (Grant No. 2022KJ174).

\end{document}